\documentclass[12pt]{amsart}
\usepackage{graphicx}
\usepackage{amssymb}
\usepackage{amsfonts}
\usepackage{amsmath}
\usepackage{array}
\usepackage{rotating}

\usepackage[all]{xy}

\newtheorem{example}{Example}[section]

\newtheorem{theorem}[example]{Theorem}

\newtheorem{conjecture}[example]{Conjecture}

\def\Proof{\noindent \it Proof -- \rm}
\def\qed{\hspace{3.5mm} \hfill \vbox{\hrule height 3pt depth 2 pt width 2mm}
\bigskip}

\def\FQSym{{\bf FQSym}}

\def\Sym{{\bf Sym}}

\def\PP{\mathring{\mathcal P}}
\def\GP{{\mathcal P}}   
\def\PQ{\mathring{\mathcal Q}}
\def\GQ{{\mathcal Q}}

\def\D{{\rm D}}
\def\<{\langle}
\def\>{\rangle}

\def\K{\operatorname{\mathbb C}}

\def\F{{\bf F}}
\def\S{{\bf S}}

\def\G{{\bf G}}
\def\M{{\bf M}}

\def\SG{{\mathfrak S}}
\def\H{{\mathcal H}}

\def\Des{\operatorname{Des}}

\def\dim{{\rm dim}}

\def\MR{{\bf MR}}
\def\bA{{\bar A}}
\def\BFQSym{{\bf BFQSym}}
\def\BSym{{\bf BSym}}

\def\maj{{\rm maj}}
\def\bmaj{{\rm bmaj}}
\def\std{{\rm std}}
\def\ol{\overline}

\title{Unital versions of the higher order peak algebras}

\author[M. Aguiar, J.-C.~Novelli,  J.-Y.~Thibon]%
{Marcelo Aguiar, Jean-Christophe Novelli, Jean-Yves Thibon}

\address[Aguiar]{Department of Mathematics\\ Texas A\&M University\\ College
Station, TX 77843\\ USA}
\address[Novelli, Thibon] {Institut Gaspard Monge, Universit\'e Paris-Est
Marne-la-Vall\'ee, \\ 5 Boulevard Descartes \\ Champs-sur-Marne \\
77454 Marne-la-Vall\'ee cedex 2 \\ France}

\email[Marcelo Aguiar]{maguiar@math.tamu.edu}
\email[Jean-Christophe Novelli]{novelli@univ-mlv.fr (corresponding author)}
\email[Jean-Yves Thibon]{jyt@univ-mlv.fr} 
\date{\today}

\begin{document}

\begin{abstract}
We construct unital extensions of the higher order peak algebras defined by
Krob and the third author in [Ann. Comb. 9 (2005), 411--430.], and show that
they can be obtained as homomorphic images of certain subalgebras of the
Mantaci-Reutenauer algebras of type $B$.
This generalizes a result of Bergeron, Nyman and the first author [Trans. AMS
356 (2004), 2781--2824.].
\end{abstract}

\maketitle

\section{Introduction}

A  \emph{descent} of a permutation $\sigma\in\SG_n$ is an index $i$
such that $\sigma(i)>\sigma(i+1)$. A descent is a \emph{peak} if moreover
$i>1$ and $\sigma(i)>\sigma(i-1)$.
The sums of permutations with a given descent set span a subalgebra of the
group algebra, the \emph{descent algebra} $\Sigma_n$.
The \emph{peak algebra} $\PP_n$ of $\SG_n$ is a subalgebra of its descent
algebra, spanned by sums of permutations having the same peak set. This
algebra has no unit.

Descent algebras can be defined for all finite Coxeter groups \cite{So}.
In \cite{ABN}, it is shown that the peak algebra of $\SG_n$ can be naturally
extended to a unital algebra, which is obtained as a homomorphic image of the
descent algebra of the hyperoctahedral group $B_n$.

The direct sum of the peak algebras turns out to be a Hopf subalgebra of the
direct sum of all descent algebras, which can itself be identified with
$\Sym$, the Hopf algebra of noncommutative symmetric functions \cite{NCSF1}.
As explained in~\cite{BHT}, it turns out that a fair amount of results on
the peak algebras can be deduced from  the case $q=-1$ of a $q$-identity
of~\cite{NCSF2}. Specializing $q$ to other roots of unity, Krob and the third
author introduced and studied \emph{higher order peak algebras} in~\cite{KT}.
Again, these are non-unital, and it is natural to ask whether the construction
of~\cite{ABN} can be extended to this case.


We will show that this is indeed possible. We first construct the unital
versions of the higher order peak algebras by a simple manipulation of
generating series. We then show that they can be obtained as homomorphic
images of the \emph{Mantaci-Reutenauer algebras} of type $B$.
Hence no Coxeter groups other than $B_n$ and $\SG_n$ are involved in the
process; in fact, the construction is related to the notion of superization,
as defined in~\cite{NT-super}, rather than to root systems or wreath products.

{\footnotesize
{\it Acknowledgements.}

This work started during a stay of the first author at the University
Paris-Est Marne-la-Vall\'ee.
Aguiar is partially supported by NSF grant DMS-0600973.
Novelli and Thibon are partially supported by the grant ANR-06-BLAN-0380.
The authors would also like to thank the contributors of the MuPAD project,
and especially of the combinat part, for providing the development environment
for their research (see~\cite{HTm} for an introduction to MuPAD-Combinat).
}

\section{Notations and background}

\subsection{Noncommutative symmetric functions}

This article is a continuation of \cite{KT}.
We will assume familiarity with the notations of \cite{NCSF1} and with the
main results of \cite{KT}. We recall a few definitions for the convenience of
the reader.

The Hopf algebra of noncommutative symmetric functions is denoted by $\Sym$,
or by $\Sym(A)$ if we consider the realization in terms of an auxiliary
alphabet $A$. Linear bases of $\Sym_n$ are labelled by compositions
$I=(i_1,\ldots,i_r)$ of $n$ (we write $I\vDash n$). 
The noncommutative complete and elementary functions are denoted
by $S_n$ and $\Lambda_n$, and $S^I=S_{i_1}\cdots S_{i_r}$.
The ribbon basis is denoted by $R_I$.
The \emph{descent set} of $I$ is
$\Des(I) = \{ i_1,\ i_1+i_2, \ldots , i_1+\dots+i_{r-1}\}$.
The \emph{descent composition} of a permutation $\sigma\in\SG_n$
is the composition $I=D(\sigma)$ of $n$ whose descent set is the descent set
of $\sigma$.

Recall from \cite{NCSF6} that for an infinite totally ordered alphabet $A$,
$\FQSym(A)$ is the subalgebra of $\K\<A\>$ spanned by the polynomials
\begin{equation}
\G_\sigma(A)=\sum_{\std(w)=\sigma}w,
\end{equation}
that is, the sum of all words in $A^n$ whose standardization is the
permutation $\sigma\in\SG_n$.
The noncommutative ribbon Schur function $R_I\in\Sym$ is then
\begin{equation}
R_I=\sum_{\D(\sigma)=I}\G_\sigma\,.
\end{equation}
This defines a Hopf embedding $\Sym\rightarrow\FQSym$.
The Hopf algebra $\FQSym$ is self-dual under the pairing 
$(\G_\sigma\,,\,\G_\tau)=\delta_{\sigma,\tau^{-1}}$
(Kronecker symbol).
Let $\F_\sigma:=\G_{\sigma^{-1}}$, so that $\{\F_\sigma\}$ is the dual basis
of $\{\G_\sigma\}$.

The {\em internal product} $*$ of $\FQSym$ is induced by
composition $\circ$ in $\SG_n$ in the basis~$\F$, that is,
\begin{equation}
\F_\sigma * \F_\tau = \F_{\sigma\circ\tau}\quad\text{and}\quad
\G_\sigma * \G_\tau = \G_{\tau\circ\sigma}\,.
\end{equation}
Each subspace $\Sym_n$ is stable under this operation, and anti-isomorphic to
the descent algebra $\Sigma_n$ of $\SG_n$.

For $f_i\in\FQSym$ and $g\in\Sym$, we have the splitting formula
\begin{equation}
\label{split}
(f_1\dots f_r)*g = \mu_r\cdot (f_1\otimes\dots\otimes f_r)*_r \Delta^r g\,,
\end{equation}
where $\mu_r$ is $r$-fold multiplication, and $\Delta^r$ the iterated
coproduct with values in the $r$-th tensor power.

\subsection{The Mantaci-Reutenauer algebra of level 2}

We denote by $\MR$ the free product $\Sym\star\Sym$ of two copies
of the Hopf algebra of noncommutative symmetric functions \cite{MR}. That is,
$\MR$ is the free associative algebra on two sequences $(S_n)$ and
$(S_{\bar n})$ ($n\ge 1$). 
We regard the two copies of $\Sym$ as noncommutative symmetric functions on 
two auxiliary alphabets: $S_n=S_n(A)$ and $S_{\bar n}=S_n(\bA)$.
We denote by $F\mapsto \bar F$ the involutive automorphism which
exchanges $S_n$ and $S_{\bar n}$. 
The bialgebra structure is defined by the requirement that the series
\begin{equation}
\sigma_1=\sum_{n\ge 0}S_n \ \text{and}\ \bar\sigma_1=\sum_{n\ge 0}S_{\bar n}
\end{equation}
are grouplike.

The internal product of $\MR$ can be computed from the splitting formula
and the conditions that $\sigma_1$ is neutral,
$\bar\sigma_1$ is central, and
$\bar\sigma_1*\bar\sigma_1=\sigma_1$.

In \cite{NT}, an embedding of $\MR$ in the Hopf algebra $\BFQSym$
of free quasi-symmetric functions of type $B$ (spanned by
colored permutations) is described. Under this embedding,
left $*$-multiplication by $\Lambda_n=\G_{n\,n-1\dots 2,1}$ corresponds to
right multiplication by $n\,n-1\dots 2,1$ in the group
algebra of $B_n$. This implies that left $*$-multiplication
by $\lambda_1$ is an involutive anti-automorphism of
$\BFQSym$, hence of $\MR$.

\subsection{Noncommutative symmetric functions of type $B$}

The hyperoctahedral analogue $\BSym$ of $\Sym$, defined in~\cite{Chow},
is the right $\Sym$-module 
freely
generated by another sequence $(\tilde S_n)$
($n\ge 0$, $\tilde S_0=1$) of homogeneous elements, with $\tilde\sigma_1$
grouplike. This is a coalgebra, but not an algebra.
It is endowed with an internal product, for which each homogeneous
component $\BSym_n$ is anti-isomorphic to the descent algebra of $B_n$.

\section{Solomon descent algebras of type $B$}

\subsection{Descents in $B_n$}

The hyperoctahedral group $B_n$ is the group of signed permutations. A signed
permutation can be denoted by $w=(\sigma,\epsilon)$ where $\sigma$ is an
ordinary permutation and $\epsilon\in\{\pm 1\}^n$, such that
$w(i)=\epsilon_i\sigma(i)$.
If we set $w(0)=0$, then, $i\in[0,n-1]$ is a descent of $w$ if $w(i)>w(i+1)$.
Hence, the descent set 
of $w$
is a subset
$D=\{i_0,i_0+i_1,\ldots,i_0+i_1+\cdots i_{r-1}\}$ of $[0,n-1]$.
We then associate to $D$ a so-called type-$B$ composition (a composition whose
first part can be zero) $(i_0-0,i_1,\dots,i_{r-1},n-i_{r-1})$.

For example, if one encodes $\epsilon$ as a boolean vector for readability,
the signed permutation $w=(231546,100100)$ has as type-$B$ composition
$I=(0,2,1,3)$.
The signed permutation $w=(231546,000100)$ has as type-$B$ composition
$I=(2,3,1)$.

The sum of all signed permutations whose descent set is contained in $D$
is mapped to $\tilde{S^I}:= \tilde S_{i_0} S^{I'}$ by Chow's
anti-isomorphism~\cite{Chow}, where $I'=(i_1,\dots,i_r)$.

\subsection{Noncommutative supersymmetric functions}

An embedding of $\BSym$ as a sub-coalgebra and sub-$\Sym$-module of $\MR$ can
be deduced from~\cite{MR}. To describe it, let us define, for
$F\in \Sym(A)$,
\begin{equation}
F^\sharp = F(A|\bar A) = F(A-q\bar A)|_{q=-1}
\end{equation}
(the supersymmetric version of $F$). 
The superization of $F\in\Sym(A)$ 
can also be given by
\begin{equation}
F^\sharp = F*\sigma_1^\sharp\,.
\end{equation}
Indeed, $\sigma_1^\sharp$ is grouplike,
and for $F=S^I$, the splitting formula gives
\begin{equation}
(S_{i_1}\cdots S_{i_r})*\sigma_1^\sharp
=\mu_r[(S_{i_1}\otimes\cdots\otimes S_{i_r})*
(\sigma_1^\sharp\otimes\cdots\otimes\sigma_1^\sharp)]=S^{I\sharp}\,.
\end{equation}

We have
\begin{equation}
\sigma_1^\sharp = \bar\lambda_1\sigma_1=\sum\Lambda_{\bar i}S_j \,.
\end{equation}
The element $\bar\sigma_1$ is central for the internal product, and
\begin{equation}
\bar\sigma_1 * F = \bar F = F*\bar\sigma_1\,.
\end{equation}
Hence,
\begin{equation}
\bar\sigma_1*\sigma_1^\sharp = \lambda_1\bar\sigma_1 =:\sigma_1^\flat\,.
\end{equation}

The basis element $\tilde S^I$ of $\BSym$, where $I=(i_0,i_1,\ldots,i_r)$ is a
type $B$-composition, can be embedded as
\begin{equation}
\tilde S^I = S_{i_0}(A)S^{i_1i_2\cdots i_r}(A|\bar A)\,.
\end{equation}
We will identify $\BSym$ with its image under this embedding.

\subsection{A proof that $\BSym$ is $*$-stable}
\label{bsym-stab}

We are now in a position to understand why $\BSym$ is a
$*$-subalgebra of $\MR$. The argument will be extended below to the case
of unital peak algebras.

Let $F,G\in \Sym$. We want to understand why
$\sigma_1 F^\sharp *\sigma_1 G^\sharp$ is in $\BSym$.
Using the splitting formula, we rewrite this as
\begin{equation}
\mu[(\sigma_1\otimes F^\sharp)*\Delta\sigma_1\Delta G^\sharp]
=\sum_{(G)}(\sigma_1G_{(1)}^\sharp)(F^\sharp * \sigma_1 G_{(2)}^\sharp).
\end{equation}
We now only have to show that each term $F^\sharp * \sigma_1 G_{(2)}^\sharp$
is in $\Sym^\sharp$. We may assume that $F=S^I$, and for any $G\in\Sym$,
\begin{equation}
S^{I\sharp}*\sigma_1 G^\sharp
=\sum_{(G)}
\mu_r[(S_{i_1}^\sharp\otimes\cdots\otimes S_{i_r}^\sharp) * 
      (\sigma_1 G_{(1)}^\sharp\otimes\cdots\otimes \sigma_1 G_{(r)}^\sharp)]
\end{equation}
so that it is sufficient to prove the property for $F=S_n$. Now,
\begin{equation}
\begin{split}
\sigma_1^\sharp * \sigma_1 G^\sharp & =
(\bar\lambda_1\sigma_1)*\sigma_1 G^\sharp \\
& = \sum_{(G)}(\bar\lambda_1 * \sigma_1 G_{(1)}^\sharp)(\sigma_1 G_{(2)}^\sharp)\\
& =\sum_{(G)}(\bar\sigma_1*\lambda_1 *\sigma_1G_{(1)}^\sharp)\cdot\sigma_1\cdot
G_{(2)}^\sharp\\
\end{split}
\end{equation}
Now,
\begin{equation}
\lambda_1 *\sigma_1G_{(1)}^\sharp
= (\lambda_1 *G_{(1)}^\sharp) (\lambda_1*\sigma_1)
= (\lambda_1 *G_{(1)}^\sharp) \lambda_1,
\end{equation}
since $\lambda_1$ is an anti-automorphism.
We then get
\begin{equation}
\begin{split}
\sigma_1^\sharp * \sigma_1 G^\sharp
& = \sum_{(G)}(\bar\sigma_1 * ((\lambda_1 *G_{(1)}^\sharp) \lambda_1)
\cdot \sigma_1\cdot G_{(2)}^\sharp\\
& =\sum_{(G)}(\bar\sigma_1*\lambda_1 *G_{(1)}^\sharp)
   \cdot(\bar\sigma_1*\lambda_1)\sigma_1\cdot
G_{(2)}^\sharp\\
& =\sum_{(G)}(\bar\lambda_1 *G_{(1)}^\sharp)\cdot\sigma_1^\sharp\cdot
G_{(2)}^\sharp\\
\end{split}
\end{equation}
Now, the result will follow if we can prove that $\bar\lambda_1 *G^\sharp$ is
in $\Sym^\sharp$ for any $G\in\Sym$.

For $G=S^I$, 
\begin{equation}
\bar\lambda_1*S^{I\sharp}
=\lambda_1*\bar\sigma_1*S^I*\sigma_1^\sharp
=\lambda_1*S^I*\bar\sigma_1*\sigma_1^\sharp
=\lambda_1*S^I*\sigma_1^\flat\,.
\end{equation}
Since left $*$-multiplication by $\lambda_1$ in an anti-automorphism,
we only need to prove that $\bar\lambda_1*S_n^{\flat}$ is of the form
$G^\sharp$. And indeed,
\begin{equation}
\begin{split}
\bar\lambda_1*S_n^{\flat}
&=\sum_{i+j=n}\lambda_1*(\Lambda_i S_{\bar j}) \\
&=\sum_{i+j=n}(\lambda_1*S_{\bar j})(\lambda_1*\Lambda_i)\\
&=\sum_{i+j=n} \Lambda_{\bar j} S_i=S_n^\sharp\,.
\end{split}
\end{equation}
This concludes the proof that $\BSym$ is a $*$-subalgebra of $\BFQSym$.

\section{Unital versions of the higher order peak algebras}

\subsection{}

As shown in \cite{BHT}, much of the theory of the peak algebra can be deduced
from a formula of \cite{NCSF2} for $R_I((1-q)A)$, in the special case $q=-1$.
In \cite{KT}, this formula was studied in the case where $q$ is an arbitrary
root of unity, and higher order analogs of the peak algebra were obtained.
In~\cite{ABN}, it was shown that the classical peak algebra can be extended to
a unital algebra, which is obtained as a homomorphic image of the descent
algebra of type $B$.

In this section, we construct unital extensions of the higher order
peak algebras.

\subsection{}
\label{s:higherpeakalg}

Let $q$ be a primitive $r$-th root of unity. All objects introduced below
will depend on $q$ (and $r$), although this dependence will not be
made explicit in the notation.

We denote by $\theta_q$ the endomorphism of $\Sym$ defined by
\begin{equation}
\tilde f =\theta_q(f)=f((1-q)A)=f(A)*\sigma_1((1-q)A)\,.
\end{equation}
We denote by $\PP$ the image of $\theta_q$ and by $\GP$ the
right $\PP$-module generated by the $S_n$ for $n\ge 0$.
Note that $\PP$ is by definition a left $*$-ideal of $\Sym$.

\begin{theorem}
$\GP$ is a unital $*$-subalgebra of $\Sym$.
Its Hilbert series is
\begin{equation}
\sum_{n\geq0}\dim \GP_n t^n = \frac{1}{1-t-t^2-\dots-t^r}\,.
\end{equation}
\end{theorem}

\Proof
Since the internal product of homogeneous elements of different degrees is
zero, it is enough to show that, for any $f,g\in\Sym$,
$\sigma_1\tilde f * \sigma_1 \tilde g$ is in $\GP$.
Thanks to the splitting formula,
\begin{equation}
\begin{split}
\sigma_1\tilde f * \sigma_1 \tilde g&=
\mu[(\sigma_1\otimes \tilde f)*\sum_{(g)}\sigma_1\tilde
g_{(1)}\otimes\sigma_1\tilde g_{(2)}]\\
&=\sum_{(g)}(\sigma_1\tilde g_{(1)})(\tilde f*\sigma_1\tilde g_{(2)})\,.
\end{split}
\end{equation}
Thus, it is enough to check that $\tilde f*\sigma_1\tilde h$ is in
$\PP$ for any $f,h\in\Sym$.
Now,
\begin{equation}
\tilde f*\sigma_1\tilde h = f*\sigma_1((1-q)A)*\sigma_1\tilde h\,,
\end{equation}
and since $\PP$ is a $\Sym$ left $*$-ideal, we only have to show that
$\sigma_1((1-q)A)*\sigma_1\tilde h$ is in $\PP$. One more splitting yields

\begin{equation}
\begin{split}
\sigma_1((1-q)A)*\sigma_1\tilde h&= (\lambda_{-q}\sigma_1)*\sigma_1\tilde h\\
&= \mu[(\lambda_{-q}\otimes\sigma_1)*\sum_{(h)}\sigma_1\tilde h_{(1)}
\otimes \sigma_1\tilde h_{(2)}]\\
&=\sum_{(h)}(\lambda_{-q}*\sigma_1\tilde h_{(1)})(\sigma_1 \tilde h_{(2)})\\
&=\sum_{(h)}(\lambda_{-q}*\tilde h_{(1)})\lambda_{-q}\sigma_1\tilde h_{(2)}
\end{split}
\end{equation}
(since left $*$-multiplication by $\lambda_{-q}$ is an anti-automorphism,
namely the composition of the antipode and $q^{\rm degree}$).
The first parentheses $(\lambda_{-q}*\tilde h_{(1)})$ are in $\PP$
since it is a left $*$-ideal. The middle term is $\sigma_1((1-q)A)$,
and the last one is in $\PP$ by definition.

Recall from \cite[Prop. 3.5]{KT} that the Hilbert series of $\PP$ is
\begin{equation}
\sum_{n\ge 0} \dim\PP_n t^n =\frac{1-t^r}{1-t-t^2-\ldots -t^r}\,.
\end{equation}
{}From \cite[Lemma 3.13 and Eq. (3.9)]{KT}, it follows that
$S_n\in\PP$ 
if and only if
$n\equiv 0\mod r$, so that the Hilbert series of $\GP$ is
\begin{equation}
\sum_{n\ge 0} \dim\GP_n t^n =\frac{1}{1-t-t^2-\ldots -t^r}\,.
\end{equation}
\qed

\section{Back to the Mantaci-Reutenauer algebra}

The above proofs are in fact special cases of a master calculation in the
Mantaci-Reutenauer algebra, which we carry out in this section.

\subsection{The $\sharp$ transform}

Let $q$ be an arbitrary complex number or an indeterminate, and define, for
any $F\in\MR$,
\begin{equation}\label{e:MRsharp}
F^\sharp = F * \sigma_1(A-q\bar A) = F * \sigma_1^\sharp\,.
\end{equation}
Since $\sigma_1^\sharp$ is grouplike, it follows from the splitting formula
that 
\begin{equation}
F\mapsto F^\sharp
\end{equation}
is an automorphism of $\MR$ for the Hopf structure. In addition, it is clear
from the definition that it is also a endomorphism of left $*$-modules.
We refer to it as the \emph{$\sharp$ transform}.

\subsection{Definition of the subalgebras}

We define
\begin{equation}
\PQ=\MR^\sharp,
\end{equation}
the image of the $\sharp$ transform. Since the latter is an endomorphism of
Hopf algebras and of left $*$-modules, $\PQ$ is both a Hopf subalgebra of
$\MR$ and a left $*$-ideal.
When $q$ is a root of unity, its image under the specialization $\bar A=A$ is
the non-unital peak algebra $\PP$ of Section~\ref{s:higherpeakalg} (and for
generic $q$, it is $\Sym$).

Let $\GQ$ be the right $\PQ$-module generated by the $S_n$, for all $n\ge 0$.
Clearly, the identification $\bar A=A$ maps $\GQ$ onto $\GP$, the unital 
peak algebra of Section~\ref{s:higherpeakalg}.

\begin{theorem}
$\GQ$ is a $*$-subalgebra of $\MR$, containing $\PQ$ as a left ideal.
\end{theorem}

\Proof
Let $F,G\in\MR$. As above, we want to show that
$\sigma_1 F^\sharp *\sigma_1 G^\sharp$ is in $\GQ$.
Using the splitting formula, we rewrite this as
\begin{equation}
\mu[(\sigma_1\otimes F^\sharp)*\Delta\sigma_1\Delta G^\sharp]
=\sum_{(G)}(\sigma_1G_{(1)}^\sharp)(F^\sharp * \sigma_1 G_{(2)}^\sharp)
\end{equation}
and we only have to show that each term $F^\sharp * \sigma_1 G_{(2)}^\sharp$
is in $\PQ$. 
We may assume that $F=S^I$, where $I$ is now a bicolored composition,
and for any $G\in\MR$,
\begin{equation}
S^{I\sharp}*\sigma_1 G^\sharp
=\sum_{(G)}
\mu_r[(S_{i_1}^\sharp\otimes\cdots\otimes S_{i_r}^\sharp) *
      (\sigma_1 G_{(1)}^\sharp\otimes\cdots\otimes \sigma_1 G_{(r)}^\sharp)]
\end{equation}
so that it is sufficient to prove the property for $F=S_n$
or $S_{\bar n}$. Now,
\begin{equation}
\begin{split}
\sigma_1^\sharp * \sigma_1 G^\sharp & =
(\bar\lambda_{-q}\sigma_1)*\sigma_1 G^\sharp \\
& = \sum_{(G)}(\bar\lambda_{-q}1 * \sigma_1 G_{(1)}^\sharp)(\sigma_1
G_{(2)}^\sharp)\\
& =\sum_{(G)}(\bar\lambda_{-q} *G_{(1)}^\sharp)\cdot\sigma_1^\sharp\cdot
G_{(2)}^\sharp
\end{split}
\end{equation}
which is in $\PQ$, since it is a subalgebra and a left $*$-ideal,
and similarly,
\begin{equation}
\begin{split}
\bar\sigma_1^\sharp * \sigma_1 G^\sharp & =
(\lambda_{-q}\bar\sigma_1)*\sigma_1 G^\sharp \\
& = \sum_{(G)}(\lambda_{-q} * \sigma_1 G_{(1)}^\sharp)(\bar\sigma_1
\bar G_{(2)}^\sharp)\\
& =\sum_{(G)}(\lambda_{-q} *G_{(1)}^\sharp)\cdot\bar\sigma_1^\sharp\cdot
\bar G_{(2)}^\sharp
\end{split}
\end{equation}
is also in $\PQ$.
\qed

The various algebras introduced in this paper and their interrelationships
are summarized in the following diagram.
\begin{equation}
\xymatrix@C=0pc{ 
{\PQ}\ar@{->>}[d]_{} &\subseteq &
{\GQ}\ar@{->>}[d]_{} &\subseteq& {\MR}\ar@{->>}[d]^{} &\subseteq&
 {\BFQSym}\ar@{->>}[d]^{}\\
 {\PP} &\subseteq& {\GP}& \subseteq & {\Sym} &\subseteq& {\FQSym}
  }
\end{equation}

Note that in the special case $q=-1$, by the results of
Section~\ref{bsym-stab}, $\GQ_n$ is the (Solomon) descent algebra of $B_n$,
$\GQ$ is isomorphic to $\BSym$, and $\GP$ is the unital peak algebra
of~\cite{ABN}.

\section{Further developments}

\subsection{Inversion of the generic $\sharp$ transform}

For generic $q$, the endomorphism~\eqref{e:MRsharp} of $\MR$ is invertible;
therefore
\begin{equation}
\PQ\sim\MR.
\end{equation}
The inverse endomorphism of $\MR$ arises from the transformation of alphabets
\begin{equation}
A\mapsto (q\bar A+A)/(1-q^2),
\end{equation}
which is to be understood in the following sense:
\begin{equation}
\label{invdiese}
\sigma_1\left(\frac{q\bar A+A}{1-q^2}\right) :=
\prod_{k\ge 0}\sigma_{q^{2k+1}}(\bar A)\sigma_{q^{2k}}(A)\,.
\end{equation}
Indeed,
\begin{equation}
\begin{split}
\sigma_1\left(\frac{q\bar A+A}{1-q^2}\right)*\sigma_1(A-q\bar A)
&= \prod_{k\ge 0}\sigma_{q^{2k+1}}(\bar A-qA)\sigma_{q^{2k}}(A-q\bar A)\\
&= \prod_{k\ge 0}\lambda_{-q^{2k+2}}(A)\sigma_{q^{2k+1}}(\bar
A)\lambda_{-q^{2k+1}}(\bar A)\sigma_{q^{2k}}(A)\\
&=\sigma_1(A)\,.
\end{split}
\end{equation}

By normalizing the term of degree $n$ in (\ref{invdiese}), we obtain
$B_n$-analogs of the $q$-Klyachko elements defined in \cite{NCSF1}:
\begin{equation}
K_n(q;A,\bar A)
:=\prod_{i=1}^n(1-q^{2\,i})S_n\left(\frac{q\bar A+A}{1-q^2}\right)
=\sum_{I\vDash n}q^{2\,\maj (I)}R_I(q\bar A+A)\,.
\end{equation}
This expression can be completely expanded on signed ribbons.
From the expression of $R_I$ in $\FQSym$, we have
\begin{equation}
R_I(\bar A+A)=\sum_{C(\sigma)=I}\G_\sigma(\bar A +A)
\end{equation}
where $\bar A +A$ is the ordinal sum. If we order $\bar A$ by
\begin{equation}
\bar a_1 < \bar a_2 < \ldots < \bar a_k < \ldots
\end{equation}
then, arguing as in \cite{NT-super}, we have 
\begin{equation}
\G_\sigma(\bar A +A)=\sum_{\std(\tau,\epsilon)=\sigma} \G_{\tau,\epsilon}
\end{equation}
so that
\begin{equation}
R_I(\bar A +A)=\sum_{\rho ({\rm J})=I} R_{\rm J}
\end{equation}
where for a signed composition ${\rm J}=(J,\epsilon)$, the unsigned
composition $\rho ({\rm J})$ is defined as the shape of
$\std(\sigma,\epsilon)$, where $\sigma$ is any permutation of shape $J$.

\subsection{}

Replacing $\bar A$ by $q\bar A$, one obtains the expansion of the
$q$-Klyachko elements of type $B$:
\begin{equation}
K_n(q;A,\bar A)=\sum_{{\rm J}}q^{\bmaj({\rm J})}R_{\rm J}
\end{equation}
where 
\begin{equation}
\bmaj({\rm J})=2\,\maj(\rho ({\rm J}))+|\epsilon|\,,
\end{equation}
where $|\epsilon|$ is the number of minus signs in $\epsilon$.

For example,
\begin{equation}
K_2(q) = R_2 + q^2\,R_{\ol2} + q^2\,R_{11} + q^3\,R_{1\ol1} +
         q\,R_{\ol11} + q^4\,R_{\ol1\ol1}.
\end{equation}

\begin{equation}
\begin{split}
K_3(q) &= R_3 + q^3\,R_{\ol3}
          + q^4\,R_{21} + q^5\,R_{2\ol1} + q^2\,R_{\ol21} + q^7\,R_{\ol2\ol1}
          + q^2\,R_{12} + q^4\,R_{1\ol2}  \\
       &  + q\,R_{\ol12} + q^5\,R_{\ol1\ol2}
          + q^6\,R_{111} + q^7\,R_{11\ol1} + q^3\,R_{1\ol11}
          + q^8\,R_{1\ol1\ol1}  \\
       &  + q^5\,R_{\ol111} + q^6\,R_{\ol11\ol1} + q^4\,R_{\ol1\ol11}
          + q^9\,R_{\ol1\ol1\ol1}.
%
%
%
\end{split}
\end{equation}

This major index of type $B$ is the flag major index defined in \cite{AR}.
Following~\cite{AR} and considering the signed
composition (where $\epsilon$ is encoded as boolean vector for readability)
\begin{equation}
{\rm J}=(2, 1, 1, \bar 3, \bar 1, \bar 2, 4, \bar 1, 2, 2)
= (2 1 1 3  1  2 4  1 2 2, 00001111110000100000)
\end{equation}
we can take the smallest permutation of shape
$(2, 1, 1, 3,  1,  2, 4,  1, 2, 2)$,
which is
\begin{equation}
\alpha=1\,5\,4\,3\,2\,6\,9\,8\,7\,11\,10\,12\,13\,16\,15\,14\,18\,17\,19
\end{equation}
sign it according to $\epsilon$, which yields
\begin{equation}
1\,5\,4\,3\,\bar 2\,\bar 6\,\bar 9\,\bar 8\,\bar
7\,\overline{11}\,10\,12\,13\,16\,\overline{15}\,14\,18\,17\,19
\end{equation}
whose standardized is
\begin{equation}
8\,11\,10\,9\,1\,2\,5\,4\,3\,6\,12\,13\,14\,16\,7\,15\,18\,17\,19
\end{equation}
and has shape $\rho({\rm J})=(2,1,1,3,1,6,3,2)$.
The major index of $\rho({\rm J})$ is $55$, the number of minus signs in
$\epsilon$ is 7, so $\bmaj({\rm J})=2\times 55+7=117$.

\subsection{}

The major index of type $B$ can be read directly on signed compositions
without reference to signed permutations as follows:
one can get $\rho({\rm J})$ by first adding the absolute values of two
consecutive parts if the left one is signed and the second one is not, then
remove the signs and proceed as before.

A different solution 
consists in reading the composition from right to left,
then associate weight $0$ (resp.~$1$) to the rightmost part if it is positive
(resp. negative) and then proceed left by adding~$2$ to the weight if the two
parts are of the same sign and~$1$ if not. Finally, add up the product of
the absolute values of the parts with their weight.

For example, with the same  ${\rm J}$ as above 
we have the following weights:
\begin{equation}
\begin{split}
{\rm J}=       & (2, 1, 1,\ \bar 3, \bar 1, \bar 2, 4, \bar 1, 2, 2) \\
{\rm weights:} &14\,12\,10\,\,  9\,\,   7\,\,\,  5\,\,4\,\,\,3\,\,\,2\,\,0
\end{split}
\end{equation}
so that we get
$2\cdot 14 + 1\cdot 12 + 1\cdot 10 + 3\cdot 9 + 1\cdot 7 + 2 \cdot 5 + 4\cdot 4 + 1 \cdot 3 + 2\cdot 2 + 2\cdot 0=117$.

This technique generalizes immediately to  colored compositions with a
fixed number $c$ of colors $0$, $1$, $\dots$, $c-1$: the weight of the
rightmost cell is its color and the weight of a part is equal to the sum of
the weight of the next part and the unique representative of the difference of
the colors of those parts modulo $c$ belonging to the interval $[1,c]$.

\subsection{Generators and Hilbert series}

For $n\ge 0$, let
\begin{equation}
S_n^\pm = S_n(A)\pm S_n(\bar A)\,,
\end{equation}
and denote by $\H_n$ the subalgebra of $\MR$ generated by the $S_k^\pm$ for
$k\le n$. For $n\ge 0$, we have
\begin{equation}
(S_n^\pm)^\sharp \equiv (1\mp q^n)S_n^\pm\ \mod\, \H_{n-1}\,,
\end{equation}
so that the $(S_n^\pm)^\sharp$ such that $1\mp q^n\not=0$ form a
set of free generators in $\MR^\sharp$.

\begin{conjecture}
If $r$ is odd, a basis of $\MR^\sharp$ will be parametrized by colored
compositions such that parts of color $0$ are not $\equiv 0 \mod r$ and parts
of color $1$ are arbitrary. The Hilbert series is then
\begin{equation}
H_r(t)=\frac{1-t^r}{1-2(t+t^2+\cdots+t^r)}.
\end{equation}
If $r$ is even, there is the extra condition that parts of color $1$ are not
$\equiv r/2 \mod r$. The Hilbert series is then
\begin{equation}
H_r(t)=\frac{1-t^r}{1-2(t+t^2+\cdots+t^r)+t^{r/2}}\,.
\end{equation}
\end{conjecture}

For example,
\begin{equation}
H_2(t)=1+t+2\,{t}^{2}+4\,{t}^{3}+8\,{t}^{4}+16\,{t}^{5}+32\,{t}^{6}+64\,{t}^{7}+128\,{t}^{8}+O
\left( {t}^{9} \right) 
\end{equation}
\begin{equation}
H_3(t)=1+2\,t+6\,{t}^{2}+17\,{t}^{3}+50\,{t}^{4}+146\,{t}^{5}+426\,{t}^{6}+1244\,{t}^{7}+3632\,{t}^{8}+O
 \left( {t}^{9} \right) 
\end{equation}
\begin{equation}
H_4(t)=1+2\,t+5\,{t}^{2}+14\,{t}^{3}+38\,{t}^{4}+104\,{t}^{5}+284\,{t}^{6}+776\,{t}^{7}+2120\,{t}^{8}+O
 \left( {t}^{9} \right)
\end{equation}

If these conjectures are correct, the Hilbert series of the
right $\MR^\sharp$-modules generated by the $S_n$ are
respectively
\begin{equation}
\frac{1}{1-2(t+t^2+\ldots+ t^r)}\,,
\end{equation}
or
\begin{equation}
\frac{1}{1-2(t+t^2+\ldots+ t^r)+t^{r/2}}\,.
\end{equation}
according to whether $r$ is odd or even.

The cases $r=1$ and $r=2$ are easily proved as follows. Assume first
that $q=1$. Set
\begin{eqnarray}
f &= 1+(\sigma_1^+)^\sharp =(\sigma_1+\lambda_{-1})(A-\bar A)\,,\\
g &= (\sigma_1^-)^\sharp-1=(\sigma_1-\lambda_{-1})(A-\bar A)\,.
\end{eqnarray}
Then, $f^2=g^2+4$, so that
\begin{equation}
f = 2\left(1+\frac14 g^2\right)^{\frac12}
\end{equation}
which proves that the $(S_n^+)^\sharp$ can be expressed
in terms of the $(S_m^-)^\sharp$.

Similarly, for $q=-1$, one can express
\begin{equation}
f = \sum_{n\ge 1}(S_{2n}^+)^\sharp +\sum_{n\ge 0}(S_{2n+1}^-)^\sharp
\end{equation}
in terms of
\begin{equation}
g = \sum_{n\ge 1}(S_{2n}^-)^\sharp +\sum_{n\ge 0}(S_{2n+1}^+)^\sharp
\end{equation}
since, as is easily verified,
\begin{equation}
(f+2)^2=g^2+4\,, \ {\it i.e., \ } f
=-2+2\left(1+\frac14 g^2\right)^{\frac12}\,.
\end{equation}
Apparently, this approach does not work anymore for higher roots
of unity.

\section{Appendix: monomial expansion of the $(1-q)$-kernel}

The results of \cite{NT-super,NCSF7} allow us to write
down a new expansion of 
$S_n((1-q)A)$, in terms of the monomial basis of \cite{AS}.
The special case $q=1$ gives back a curious expression of
Dynkin's idempotent, first obtained in \cite{AL}.

Let $\sigma$ be a permutation.
We then define its \emph{left-right minima} set ${\rm LR}(\sigma)$ as the
values of $\sigma$ that have no smaller value to their left.
We will denote by ${\rm lr}(\sigma)$ the cardinality of ${\rm LR}(\sigma)$.
For example, with $\sigma=46735182$, we have
${\rm LR}(\sigma)=\{4,3,1\}$, and
${\rm lr}(\sigma)=3$.

Let us now compute how $S_n((1-q)A)$ decomposes on the monomial
basis $\M_\sigma$ (see~\cite{AS}) of $\FQSym$.
Thanks to the Cauchy formula of $\FQSym$~\cite{NCSF7}, we have
\begin{equation}
S_n((1-q)A) = \sum_\sigma \S^{\sigma}(1-q) \M_\sigma(A),
\end{equation}
where $\S$ is the dual basis of $\M$. Given the transition matrix between $\M$
and $\G$, we immediately deduce that
\begin{equation}
\S^\sigma = \sum_{\tau < \sigma^{-1}} \F_\tau,
\end{equation}
where $<$ stands for the right weak order on permutations, so that, for
example.
\begin{equation}
\S^{312} = \F_{123} + \F_{213} + \F_{231}.
\end{equation}

Thanks to~\cite{NT-super}, we know that $\F_\sigma(1-q)$ is either $(-q)^k$ if
$\Des(\sigma)=\{1,\dots,k\}$ or $0$ otherwise.
Let us define \emph{hook permutations} of hook $k$ the permutations $\sigma$
such that $\Des(\sigma)=\{1,\dots,k\}$.
Now, $\S^{\sigma}(1-q)$ amounts to compute the list of \emph{hook permutations}
smaller than $\sigma$.
Note that hook permutations are completely characterized by their left-right
minima. Moreover, if $\tau$ is smaller than $\sigma$ in the right weak order,
then ${\rm LR}(\tau)\subset {\rm LR}(\sigma)$.

Hence all hook permutations smaller than a given permutation $\sigma$ belong
to the set of hook permutations with left-right minima in ${\rm LR}(\sigma)$.
Since by elementary transpositions decreasing the length, one can get from
$\sigma$ to the hook permutation with the same left-right minima and then from
this permutation to all the others, we have:

\begin{theorem}
Let $n$ be an integer.
Then
\begin{equation}
S_n((1-q)A) = \sum_{\sigma\in\SG_n} (1-q)^{{\rm lr}(\sigma)} \M_\sigma.
\end{equation}
\qed
\end{theorem}

In the particular case $q=1$, we recover a result of~\cite{AL}:
\begin{equation}
\Psi_n= \sum_{\genfrac{}{}{0pt}{}{\sigma\in\SG_n}{\sigma(1)=1}}\M_\sigma,
\end{equation}
where $\Psi_n$ is the {\em first Eulerian idempotent}~\cite[Prop. 5.2]{NCSF2}.


\end{document}